\documentclass[10pt,twoside]{amsart}
\usepackage{amsmath}
\usepackage{amssymb}
\usepackage{graphicx}
\usepackage[colorlinks=true, urlcolor=blue,bookmarks=true,bookmarksopen=true, citecolor=blue]{hyperref}

\numberwithin{equation}{section}

\newtheorem{theorem}{Theorem}[section]

\newtheorem{corollary}[theorem]{Corollary}

\newtheorem{proposition}[theorem]{Proposition}
\newtheorem{remark}[theorem]{Remark}
\def \bb{\mathbb}

\newcommand{\comb}[2]{{\left(\begin{smallmatrix} #1 \\ #2 \end{smallmatrix}\right)}}

\def \CP{{\bb{CP}}}

\def \PP{{\bb{P}}}

\def \RR{{\bb{R}}}

\def \ZZ{{\bb{Z}}}

\def \({\left(}
\def \){\right)}
\def \<{\langle}
\def \>{\rangle}

\def \bar{\overline}
\def \deg{\mathrm{deg}}
\def \dsum{\oplus}

\def \into{\hookrightarrow}

\def \tensor{\otimes}

\def \vargeq{\geqslant}

\def \Ham{{\rm Ham}}

\def \img{{\rm img }}

\def \Span{{\rm Span}}

\newcommand{\1}{1\!\!\!1}

\def \qed{\hfill $\square$ \vspace{0.03in}}
\begin{document}

\title{An example concerning Hamiltonian groups of self product, I}

\author{Shengda Hu}
\address{Department of Mathematics, Wilfrid Laurier University, 75 University Ave. West, Waterloo, Canada}
\email{shu@wlu.ca}

\author{Fran\c{c}ois Lalonde}
\address{D\'epartement de math\'ematiques et de Statistique, Universit\'e de Montr\'eal, C.P. 6128, Succ. Centre-ville, Montr\'eal H3C 3J7, Qu\'ebec, Canada}
\email{lalonde@dms.umontreal.ca}

% \begin{abstract}
\abstract
We show that $(S^2\times S^2, \omega_0 \dsum \lambda\omega_0)$, with $\lambda > 1$, is an example of symplectic manifold $(X, \omega)$ such that the $\pi_1\Ham(X \times X, \omega\dsum -\omega)$ contains extra elements than those from $\pi_1\Ham(X, \omega) \times \pi_1\Ham(X, -\omega)$. 
\endabstract
% \end{abstract}

\maketitle


\vspace{.08in} \noindent \textbf{Keywords}: Hamiltonian group, Seidel elements.

\section{Introduction}

Let $(X, \omega)$ be a compact symplectic manifold with $\dim_\RR X = 2n$ and $\Ham(X, \omega)$ the group of Hamiltonian diffeomorphisms.
It's natural to ask how $\Ham(X, \omega) \times \Ham(X, -\omega)$ compares with $\Ham(X \times X, \omega \dsum - \omega)$. Firstly, there is a natural injection:
$$m : \Ham(X, \omega) \times \Ham(X, -\omega) \into \Ham(X \times X, \omega \dsum - \omega) : m(\phi, \psi) = (\phi, \psi)$$
Secondly, since a neighbourhood of the diagonal $\triangle \subset X \times X$ is symplectomorphic to a neighbourhood of the zero section in $T^*X$, it is clear that the injection $m$ can't be surjective. On the other hand, it is not as clear how they compare homotopically. It is well known that for $(X, \omega) = (S^2, \omega_0)$, the standard $2$-sphere, the two sides of $m$ are weakly homotopic. In this article we consider the first homotopy group, and will use $m$ to denote the induced map on $\pi_1$ as well.
To save notations, we use $X$ to denote $(X, \omega)$ and $\bar X$ to denote $(X, -\omega)$.

Seidel constructed for each $\gamma \in \pi_1\Ham(X)$ an automorphism $\Phi^X_\gamma$ of the quantum homology ring $QH_*(X)$ as a module over itself. Let $\1 = [X] \in QH_*(X)$ be the unit, then the \emph{Seidel element} $\Psi^X_\gamma = \Phi^X_\gamma(\1) \in QH_*^\times(X)$ is an invertible element.\footnote{Seidel's original construction \cite{Seidel} gives for each choice of a reference section an automorphism as well as an element. Here, we follow McDuff \cite{McDuff}, choosing a canonical reference section and refer to the result as \emph{the} Seidel morphism and element. Both will appear in the main text.} The map $\Psi^X : \pi_1\Ham(X) \to QH_*^\times(X) : \Psi^X(\gamma) = \Psi^X_\gamma$ is the \emph{Seidel homomorphism}, where $QH_*^\times$ is a group under quantum multiplication. 
 
In this article, we consider the example $(X, \omega)= (S^2 \times S^2, \omega_0 \dsum \lambda\omega_0)$, where $\omega_0$ is the standard volume form on $S^2$ and $\lambda > 1$. We prove the following statement, using explicit computation of the Seidel elements.
\begin{theorem}\label{thm:main}
$m$ is not surjective on $\pi_1$ for $(X, \omega) = (S^2 \times S^2, \omega_0 \dsum \lambda\omega_0)$ with $\lambda > 1$.
\end{theorem}

\begin{remark}\label{diagonal:almostlast}
\rm{
We note that, in fact, $\pi_1\Ham(X, \omega)$ already has an element $S$ which does not come from $\pi_1\Ham$ of either of its factors. On the other hand, the factors are not symplectomorphic (after reversing one of the structures). Indeed, Gromov \cite{Gromov} showed that $\Ham(X, \omega_0 \dsum \omega_0)$ is weakly homotopic to $SO(3) \times SO(3)$, which in turn is weakly homotopic to $\Ham(S^2, \omega_0) \times \Ham(S^2, \omega_0)$.
}
\end{remark}

Let's start by fixing some notations. Let $\Gamma_\omega = \pi_2(M) / \sim$ where $\beta \sim \beta' \iff \omega(\beta - \beta') = c_1(TX)(\beta - \beta') = 0$. As a group, the quantum homology $QH_*(X, \omega) \cong H_*(X, \omega) \tensor \Lambda_\omega$ where $\Lambda_\omega$ is the Novikov ring
$$\Lambda_\omega = \left\{\left.\sum_{\beta \in \Gamma_\omega} a_\beta e^\beta\right| a_\beta \in \RR, \forall K \in \RR, \# \{\beta | a_\beta \neq 0 \text{ and } \omega(\beta) > K\} < \infty\right\}$$
graded by $\deg e^\beta = 2c_1(TX)(\beta)$.
The quantum (intersection) product on $QH_*(X)$ is given by
$$a * b = \sum_{\beta\in\Gamma_\omega, c\in H_*(X, \omega)} \<a, b, \hat c\>_\beta e^{-\beta} c$$
where $\hat c \in H_*(X)$ is the Poincar\'e dual of $c$ under the ordinary intersection product and $\<a, b, \hat c\>_\beta$ is the genus $0$ Gromov-Witten invariant counting the number of $J$-holomorphic rational curves in $X$  passing through representatives of the classes $a$, $b$ and $\hat c$, representing the class $\beta$.

Next recall the effect of reversing the symplectic structure on $QH_*(X)$ and the Seidel elements. It leaves $\Gamma_\omega$ unchanged. Let $\tau : \pi_2(X) \to \pi_2(X) : \beta \mapsto - \beta$, it induces the ring isomorphism
$$\tau : \Lambda_\omega \to \Lambda_{-\omega} : \sum_{\beta \in \Gamma_\omega} a_\beta e^\beta \mapsto \sum_{\beta \in \Gamma_{\omega}} a_\beta \tau\(e^{\beta}\) = \sum_{\beta \in \Gamma_{\omega}} (-1)^{c_1(TX)(\beta)}a_\beta e^{ -\beta}$$
The quantum homology $QH_*(X)$ and $QH_*(\bar X)$ are isomorphic as rings via
$$\tau : QH_*(X) \to QH_*(\bar X) : \tau(a \tensor e^\beta) = (-1)^{n+c_1(TX)(\beta)} a \tensor e^{-\beta}$$
where $a \in H_*(X)$.
Let $\gamma = [g] \in \pi_1\Ham(X)$ where $g \in \Omega_0\Ham(X,\omega)$ is a loop in $\Ham(X)$ based at $id$ and define $\tau : \pi_1\Ham(X) \to \pi_1\Ham(\bar X)$ by $\tau(\gamma) = [g^-]$, where $g^-(t) = g(1-t)$, then the Seidel elements are related by
\begin{equation}
\label{SeidelReversion}
\tau(\Psi^X_\gamma) = \Psi^{\bar X}_{\tau(\gamma)}
\end{equation}

Let $(X, \omega_X)$ and $(Y, \omega_Y)$ be compact monotone symplectic manifolds, then we have the ring isomorphism extending the K\"unneth isomorphism for ordinary homology:
\begin{equation}\label{eq:quantumKuneth}
QH_*(X \times Y, \omega_X \dsum \omega_Y) \cong QH_*(X, \omega_X) \tensor QH_*(Y, \omega_Y)
\end{equation}
For the case under consideration, although $(X, \omega) = (S^2 \times S^2, \omega_0 \dsum \lambda\omega_0)$ is not monotone, neither is $(X, -\omega)$, the manifold $(X \times X, \omega \dsum -\omega)$ can be written as a product of monotone manifolds:
$$(X \times X, \omega \dsum -\omega) = (X_1 \times X_1, \omega_1 \dsum \lambda \omega_1)$$
where $\omega_1 = \omega_0 \dsum -\omega_0$ on $X_1 = S^2 \times S^2$. 
Since
$$QH_*(X_1, \omega_1) \tensor QH_*(X_1, \lambda \omega_1) \cong QH_*(S^2, \omega_0) \tensor QH_*(S^2, -\omega_0) \tensor QH_*(S^2, \lambda\omega_0) \tensor QH_*(S^2, -\lambda\omega_0)$$
it follows still that
$$QH_*(X \times X, \omega \dsum -\omega) \cong QH_*(X, \omega) \tensor QH_*(X, -\omega)$$

The Hamiltonian groups are similarly related:
$$m: \Ham(X, \omega_X) \times \Ham(Y, \omega_Y) \into \Ham(X \times Y, \omega_X \dsum \omega_Y)$$
Moreover, let $\gamma_X \in \pi_1\Ham(X, \omega_X)$ and $\gamma_Y\in \pi_1\Ham(Y, \omega_Y)$ then $\gamma_{X\times Y}:= m(\gamma_X, \gamma_Y) \in \pi_1\Ham(X \times Y,\omega_X \dsum \omega_Y)$. Suppose that the ring isomorphism \eqref{eq:quantumKuneth} holds, then the respective Seidel elements are related by
$$\Psi^{X\times Y} (\gamma_{X\times Y}) = \Psi^X (\gamma_X) \tensor \Psi^Y (\gamma_Y)$$

\vspace{0.1in}\noindent
{\bf Acknowledgement.} S. Hu would like to thank Dr. S. Lisi for helpful discussions. S. Hu is partially supported by an NSERC Discovery Grant. We acknowledge the referee's suggestion to divide our paper into two parts.

\section{Example: $(X, \omega) = (S^2 \times S^2, \omega_0 \dsum \lambda \omega_0)$}\label{example}

Let $(S^2, \omega_0)$ be the sphere with the standard symplectic structure, $X = (S^2 \times S^2, \omega_0 \dsum \lambda\omega_0)$  for some $\lambda > 0$, and $(M, \Omega) = X \times \bar X$.
Denote the factors as $\PP_j$ for $j = 1, \ldots, 4$. Let 
$$(X', \omega') = \PP_1\times \PP_4 \text{ and } (M', \Omega') = X' \times \bar X',$$
then $M'$ and $M$ are isomorphic symplectic manifolds, by switching the factors; while $X'$ and $X$ are isomorphic via an anti-symplectic involution on the second factor. 

When $\lambda \in (1, 2]$, it's known (see for example McDuff-Tolman \cite{McDuffTolman}) that $\pi_1\Ham(X)$ is generated by $3$ elements: $r_1$ and $r_2$ of order $2$ rotating the respective factors and an element $s$ of infinite degree. $X$ admits another structure of $S^2$ fibration over $S^2$ and $s$ defines an $S^1$action on $X$ rotating the fibers. The diagonal and the anti-diagonal are the two sections of the fibration fixed by this $S^1$-action, and the weight of the action on the normal bundle of the section with bigger area is $-1$.

In order to write down the Seidel elements in $QH_*(X)$ and for later convenience, we introduce a system of notations for the elements in $H_*$ of the various spaces involved. The homology $H_*(S^2) = \ZZ \dsum 0 \dsum \ZZ$, as graded by the degree. We write $(1) \in H_2(S^2)$ and $(0) \in H_0(S^2)$ as the respective (positive) generators (with respect to the volume form $\omega_0$). For a (positive) basis of $H_*(S^2)$ with respect to the reverse form $-\omega_0$, we write $(\bar 1) := -(1) \in H_2(S^2)$ and $(\bar 0) := - (0) \in H_0(S^2)$.
The homology $H_*(X)$ is then generated by $(11) \in H_4(X)$, $(10), (01) \in H_2(X)$ and $(00) \in H_0(X)$, where, for example, $(10)$ denotes the tensor $(1) \tensor (0)$. We use similar notations for the generators of $H_*(M)$, e.g. $(01\bar 0\bar 1) \in H_4(M)$. 

The quantum homology $QH_*(S^2)$ is determined by the fact that $(1)$ is the unit and
$$(0) * (0) = (1) e^{-(1)}$$
For $QH_*(\bar{S^2})$, we have the corresponding $\bar{\left.\cdot \right.}$-version:
$$(\bar 0) \bar{*} (\bar 0) = (\bar 1) e^{-(\bar 1)} \Rightarrow (0) \bar{*} (0) = -(1) e^{(1)}$$
Note that the unit in the quantum homology $QH_*(X)$, $QH_*(X')$ and $QH_*(M)$ are respectively $(11)$, $(1\bar 1)$ and $(11\bar 1 \bar 1)$.
We have for example
$$(01) * (10) = (00) \text{ and } (01\bar 0 \bar 1) * (00\bar 1\bar 1) = (10\bar 0\bar 1) e^{-(1000)}$$
Using these notations, let $r$ denote the action of $S^1$ on $S^2$ fixing the poles and $\Psi_{r} \in QH_*(S^2)$ be the corresponding Seidel element, then
$$\Psi^{S^2}_{r} = (0) e^{\frac{1}{2}(1)} \text{ and } \Psi^{\bar S^2}_{\tau(r)} = \tau(\Psi^{S^2}_r) = (-1)^{c_1(TS^2)\(\frac{1}{2}(1)\)}(\bar 0) e^{-\frac{1}{2}(1)} = -(\bar 0) e^{-\frac{1}{2}(1)} \in QH_*(\bar{S^2})$$

We write down the Seidel elements for $R_1$ and $R_2$:
$$\Psi^X_{r_1} = \Psi^{S^2}_{r} \tensor \Psi^{S^2}_{\1} = (01)e^{\frac{1}{2}(10)} \text{ and } \Psi^X_{r_2} = \Psi^{S^2}_{\1} \tensor \Psi^{S^2}_{r} = (10) e^{\frac{1}{2}(01)}$$
Following \cite{McDuffTolman}, we explicitly write down the Seidel element for $s$:
$$\Psi^X_s = [(01) + (10)] e^{\frac{1}{2}(10) + h[(10) - (01)]} \text{ where } h = \frac{1}{6\lambda(\lambda - 1)}$$
where $\omega((10)) = 1$, $\omega((01)) = \lambda$ and  $c_1((01)) = c_1((10)) = 2$. Because
$$[(01) + (10)] * [(01) - (10)] = (11)\(e^{-(10)} - e^{-(01)}\)$$
we see that the reversed loop $s^-$ gives the Seidel element
$$\Psi^X_{s^-} = (\Psi^X_s)^{-1} = [(01) - (10)] e^{\frac{1}{2}(10)- h[(10) - (01)]}\(1+e^{(10)-(01)} + e^{2[(10)-(01)]} + \ldots \)$$
The corresponding Seidel elements in $QH_*(\bar X)$ are:
$$\Psi^{\bar X}_{\tau(r_1)} = -(\bar 0\bar 1)e^{-\frac{1}{2}(10)}, \Psi^{\bar X}_{\tau(r_2)} = -(\bar 1\bar 0)e^{-\frac{1}{2}(01)} \text{ and }$$
$$\Psi^{\bar X}_{\tau(s)} = -[(\bar 0\bar 1) + (\bar 1\bar 0)] e^{-\frac{1}{2}(10) - h[(10) - (01)]}.$$

Next we describe the Seidel elements in $QH_*(X')$. Those for $r_1'$ and $r_2'$ are:
$$\Psi^{X'}_{r_1'} = \Psi^{S^2}_r \tensor \Psi^{\bar{S^2}}_{\tau(\1)} = (0\bar 1) e^{\frac{1}{2}(10)} \text{ and } \Psi^{X'}_{r_2'} = \Psi^{S^2}_{\1} \tensor \Psi^{\bar{S^2}}_{\tau(r)} = -(1\bar 0)e^{-\frac{1}{2}(01)}.$$
To describe the Seidel elements of infinite order, we notice that $(X', \omega')$ is symplectically identified with $(X, \omega)$ by
$$(1, c): \CP^1 \times \CP^1 \to \CP^1 \times \CP^1,$$
where $c$ is the antipodal map. It induces on $H_*$ the isomorphism given by
$$(1, c)_*: ((00), (01), (10), (11)) \mapsto ((00), (0\bar 1), (10), (1\bar 1))$$
from which can be recovered the expressions for $\Psi^{X'}_{r_1'}$ and $\Psi^{X'}_{r_2'}$ given above.
Let $s'$ be the loop conjugate to $s$ by the map $(1,c)$ then the corresponding Seidel element is
$$\Psi^{X'}_{s'} = [(0\bar 1) - (1\bar 0)] e^{\frac{1}{2}(10) + h[(01) + (10)]} \in QH_*(X', \omega').$$
The corresponding Seidel elements in $QH_*(\bar X')$ are:
$$\Psi^{\bar X'}_{\tau(r_1')} = -(\bar 0 1) e^{-\frac{1}{2}(10)}, \Psi^{\bar X'}_{\tau(r_2')} = (\bar 1 0) e^{\frac{1}{2}(01)} \text{ and }$$
$$\Psi^{\bar X'}_{\tau(r')} = -[(\bar 0 1)-(\bar 1 0)] e^{-\frac{1}{2}(10) - h[(01) + (10)]}.$$

The image of the obvious map:
$$m: \pi_1\Ham(X) \times \pi_1\Ham(\bar X) \to \pi_1\Ham(M)$$
is generated by the image of $\{\1, r_1, r_2, s\} \times \{\1, \tau(r_1), \tau(r_2), \tau(s)\}$ and the corresponding Seidel elements are given by the respective tensor products. %, by the lemma \ref{diagonal:seploop}. 
Let $m'$ be the corresponding map for $(X', \pm\omega')$:
$$m' : \pi_1\Ham(X') \times \pi_1\Ham(\bar X') \to \pi_1\Ham(M') = \pi_1\Ham(M),$$
where the last identification is by switching the factors of $M'$.
The image of $m'$ is generated by the image of $\{\1, r_1', r_2', s'\} \times \{\1, \tau(r_1'), \tau(r_2'), \tau(s')\}$. Simple algebraic observation together with the explicit description of the Seidel elements given above lead to
\begin{proposition}\label{diagonal:extraloop}
$\img (m) \neq \img (m') \subset \pi_1\Ham(M, \Omega)$.%, in particular, $m$ is not surjective.
\end{proposition}
{\it Proof:} We first proceed as far as possible without using the exact form of the Seidel elements computed above.
Let $S = m(s, \1)$, $T = m(\1, \tau(s))$, $R_j = m(r_j, \1)$, $\bar R_j = m(\1, \tau(r_j))$ for $j = 1,2$ and the corresponding ones with $\left.\right.'$, be loops in $\Ham(M, \Omega)$. Let $\Lambda := \Lambda_\Omega$ denote the Novikov ring for $(M, \Omega)$.
It's evident that
\begin{equation}\label{diagonal:seidelform}
\begin{split}
\Psi^M_{ S} \in \Span_\Lambda((01\bar 1\bar 1), (10\bar 1\bar 1)), \hspace{0.1in} & \Psi^M_{ T} \in \Span_\Lambda((11\bar 0\bar 1), (11\bar 1\bar 0)), \text{ and }\\
\Psi^M_{ S'} \in \Span_\Lambda((01\bar 1\bar 1), (11\bar 1\bar 0)), \hspace{0.1in} & \Psi^M_{ T'}\in \Span_\Lambda((10\bar 1\bar 1), (11\bar 0\bar 1)).
\end{split}\end{equation}
More explicitly, we have the following
\begin{equation*}
\begin{split}
\Psi^M_{ S} & = \left[(01\bar 1 \bar 1) + (10\bar 1 \bar 1)\right]e^{\frac{1}{2}(1000) + h\left[(1000) - (0100)\right]}\\
\Psi^M_{ T} & = -\left[(11\bar 0 \bar 1) + (11\bar 1 \bar 0)\right]e^{-\frac{1}{2}(0010) - h\left[(0010) - (0001)\right]}\\
\Psi^M_{ S'} & = \left[-(11\bar 1 \bar 0) + (01\bar 1 \bar 1)\right]e^{\frac{1}{2}(1000) + h\left[(0001) + (1000)\right]}\\
\Psi^M_{ T'} & = -\left[-(10\bar 1 \bar 1) + (11\bar 0 \bar 1)\right]e^{-\frac{1}{2}(0010) - h\left[(0100) + (0010)\right]}\\
\end{split}
\end{equation*}

We'll drop the superscripts such as ${}^X$ from the notation of the Seidel elements as they can be inferred from the subscripts. The Seidel elements of loops in $\img(m)$ are of the form 
$$\sigma := \Psi_{ R_1}^{\epsilon_1}\Psi_{ R_2}^{\epsilon_2}\Psi_{\bar R_1}^{\epsilon_3} \Psi_{\bar R_2}^{\epsilon_4} \Psi_{ S}^{p} \Psi_{ T}^{q}$$
where $\epsilon_j \in \{0, 1\}$ and $p, q \in \ZZ$. Square it we have
\begin{equation}\label{diagonal:squareseidel}
\sigma^2 = \Psi_{ S}^{2p} \Psi_{ T}^{2q}
\end{equation}
Suppose that $\sigma$ also lies in $\img(m')$, then $\exists p', q' \in \ZZ$ so that
\begin{equation}\label{diagonal:contradict}\sigma^2 = \Psi_{ S}^{2p} \Psi_{ T}^{2q} = \Psi_{ S'}^{2p'} \Psi_{ T'}^{2q'} = {\sigma'}^2
\end{equation}
In the following we show that \eqref{diagonal:contradict} holds iff $p = q = p' = q' = 0$. 

It's easy to see from \eqref{diagonal:seidelform} (also see below for the first two) that 
\begin{equation*}
\begin{split}
\Psi_{ S}^2 \in V:= \Span_\Lambda((11\bar 1\bar 1), (00\bar 1\bar 1)), \hspace{0.1in} & \Psi_{ T}^2 \in W := \Span_\Lambda((11\bar 1\bar 1), (11\bar 0\bar 0)) \\
\text{ and } \Psi_{ S'}^2 \in V' := \Span_\Lambda((11\bar 1\bar 1), (01\bar 1\bar 0)), \hspace{0.1in} & \Psi_{ T'}^2\in W' := \Span_\Lambda((11\bar 1\bar 1), (10\bar 0\bar 1)).
\end{split}
\end{equation*}
Notice that $V, V', W$ and $W'$ are closed under the quantum product $*$ and inverse (whenever exists).

Let us first assume that $p, q, p', q' \vargeq 0$, then $\sigma^2$ has the form:
$$(a(11\bar 1\bar 1) + b(00\bar 1\bar 1))*(c(11\bar 1\bar 1) + d(11\bar 0\bar 0)) = ac(11\bar 1\bar 1) + ad(11\bar 0\bar 0) + bc(00\bar 1\bar 1) + bd(00\bar 0\bar 0)$$
while ${\sigma'}^2$ is of the form:
$$(a'(11\bar 1\bar 1) + b'(10\bar 0\bar 1)) *(c'(11\bar 1\bar 1) + d'(01\bar 1\bar 0)) = a'c'(11\bar 1\bar 1) + a'd'(01\bar 1\bar 0) + b'c'(10\bar 0\bar 1) + b'd'(00\bar 0\bar 0)$$
It follows that the necessary condition for \eqref{diagonal:contradict} to hold is
\begin{equation}\label{diagonal:coefvanish}ad = bc = a'd' = b'c' = 0 \in \Lambda
\end{equation}

Here we need the explicit form of the Seidel elements. First we have
$$\Psi_s^2 = \left[2(00) + (11) \left(e^{-(10)} + e^{-(01)}\right)\right] e^{(10) + 2h[(10) - (01)]} \in QH_*(X).$$
Now let $x = e^{-(10)}$, $y = e^{-(01)}$, $A = (00)$ and $B = (11)$, then for any integer $p > 0$
$$\Psi_s^{2p} = K^p \left(A + \frac{x+y}{2} B\right)^p, \text{ where } A^2 = Bxy, B^2 = B, AB = A
\text{ and } K = 2x^{-2h-1} y^{2h}$$ We have the explicit formula
$$\Psi_s^{2p} = K^p \left(\sum_{i = 0}^{\lfloor\frac{p}{2}\rfloor}\comb{p}{2i} \alpha^{p-2i} (xy)^i B + \sum_{i = 0}^{\lfloor\frac{p-1}{2}\rfloor} \comb{p}{2i+1} \alpha^{p-2i-1} (xy)^i A\right), \text{ where }\alpha = \frac{x+y}{2}.$$
Note that
$$\tau(x) = e^{(10)} = x^{-1}, \tau(y) = e^{(01)} = y^{-1}, \tau(A) = (\bar 0\bar 0) = (00) = A \text{ and } \tau(B) = B$$
It follows that $\tau(\alpha) = (xy)^{-1}\alpha$ and $\tau(K) = 2x^{2h+1}y^{-2h} = 4K^{-1}$. Using \eqref{SeidelReversion} we get for $q > 0$
$$\Psi_{\tau(s)}^{2q} = 4^qK^{-q}\left(\sum_{i = 0}^{\lfloor\frac{q}{2}\rfloor}\comb{q}{2i} \alpha^{q-2i} (xy)^{i-q} B + \sum_{i = 0}^{\lfloor\frac{q-1}{2}\rfloor} \comb{q}{2i+1} \alpha^{q-2i-1} (xy)^{i+1-q} A\right),$$
Since $\Psi_{ S}^{2p} = \Psi_{s}^{2p} \tensor \Psi_{\tau(\1)}$ and $\Psi_{ T}^{2q} = \Psi_{\1} \tensor \Psi_{\tau(s)}^{2q}$, it follows that in \eqref{diagonal:coefvanish} $ad = bc = 0 \Rightarrow p = q = 0$, i.e. $\sigma^2 = id$. Similaly $a'd' = b'c' = 0 \Rightarrow p' = q' = 0$ and $(\sigma')^2 = id$.

The other cases of the sign combinations of $p, q, p'$ and $q'$ are similar. Among $p, q, -p', -q'$, there must be $2$ of the same sign. Let's suppose $p$ and $-p'$ are of the same sign, say both $\vargeq 0$, then instead of \eqref{diagonal:contradict} we may consider
$$\Psi_{ S}^{2p} \Psi_{ S'}^{-2p'} = \Psi_{ T}^{-2q} \Psi_{ T'}^{2q'}.$$
Without using the details of the Seidel elements involved, we arrive at an equation similar to \eqref{diagonal:coefvanish}. Afterwards, explicit computation similar to the above gives $p = p' = 0$ and thus $\sigma^2 = (\sigma')^2 = id$.

It follows that, at least, all elements in the image of $m$ of the form $p S + q  T$ with $p$ or $q \neq 0$  do not lie in the image of $m'$, and the proposition follows.
\qed

\begin{corollary}\label{coro:mainthm}
$m$ is not surjective on $\pi_1$ for $(X, \omega) = (S^2 \times S^2, \omega_0 \dsum \lambda\omega_0)$ with $\lambda > 1$. \qed
\end{corollary}

\label{lastpage-01}
\end{document}